\documentclass[11pt,a4paper]{amsart}
\usepackage{amssymb}
\usepackage{amsmath}
\usepackage{marginnote}
\newcommand{\R}{\mathbb{R}}

\font\eufm=eufm10
\def\frak#1{\hbox{\eufm#1}}
\newcommand{\bd}{\begin{document}}
\newcommand{\ed}{\end{document}}
\newcommand{\be}{\begin{enumerate}}
\newcommand{\ee}{\end{enumerate}}
\newcommand{\bi}{\begin{itemize}}
\newcommand{\ei}{\end{itemize}}
\newcommand{\ba}{\begin{array}}
\newcommand{\ea}{\end{array}}
\newcommand{\vs}{\vspace*{0.3\baselineskip}}%%%maly odstep pionowy
\newcommand{\vsm}{\vspace*{-0.3\baselineskip}}%%%maly odstep pionowy ujemny
\newcommand{\kom}[1]{{\em #1}\newline}%%%komentarz w oddzielnej lini
\newtheorem{defi}{Definition}[section]
\newtheorem{tw}[defi]{Theorem}
\newtheorem{prop}[defi]{Proposition}
\newtheorem{lem}[defi]{Lemma}
\newtheorem{re}[defi]{Remark}
\newtheorem{col}[defi]{Corollary}
\newtheorem{ex}[defi]{Examples}
\newtheorem{zad}{Exercise}[section]
\newtheorem{zal}{Assumptions}[section]
\newtheorem{assumpt}[defi]{Assumptions}
\newcommand{\Om}{\Omega}
\newcommand{\om}{\omega}
\newcommand{\G}{\Gamma}
\newcommand{\D}{\Delta}
\renewcommand{\d}{\delta}
\newcommand{\ga}{\gamma}
\newcommand{\eps}{\epsilon}
\newcommand{\ove}{\overline}
\newcommand{\ms}{\oplus}
\newcommand{\mt}{\otimes}
\newcommand{\dz}{\wedge}
\newcommand{\lra}{\longrightarrow}
\newcommand{\sign}{\mbox{$ sgn $}}
\newcommand{\rel}{\mbox{$\,$\rule[0.5ex]{1.1em}{0.2pt}$\triangleright\,$}}
\newcommand{\dow}{\hspace*{\fill}\rule{1.6ex}{1.6ex}\hspace*{1em}}
\newcommand{\dowl}{\hspace*{\fill}\rule{1ex}{1ex}\hspace*{1em}}
\newcommand{\sd}{\hspace{0.3ex}\tiny{\rhd\mbox{\hspace{-2ex}}<}\hspace{0.3ex}}
\newcommand{\mmt}[2]{\mbox{$\vphantom{}_{#1}\times_{#2}$}}
\newcommand{\gotg}{\frak g}
\newcommand{\gota}{\frak a}
\newcommand{\gotb}{\frak b}
\newcommand{\gotc}{\frak c}
\newcommand{\gothh}{\frak h}
\newcommand{\gott}{\frak t}
\newcommand{\hd}{\hat{\d}}
\newcommand{\oml}{\Omega_L^{1/2}}
\newcommand{\omr}{\Omega_R^{1/2}}
\newcommand{\omh}{\Omega^{1/2}}
\newcommand{\lo}{\lambda_0}
\newcommand{\ro}{\rho_0}
\newcommand{\lma}{\Lambda^{max}}
\newcommand{\timh}{\times_h}
\newcommand{\Gd}{\G^{(2)}}
\newcommand{\el}{e_L}
\newcommand{\er}{e_R}
\newcommand{\GG}{\G_1\times\G_2}
\newcommand{\gdot}{\hspace{-0.1em}\cdot\hspace{-0.1em}}
\newcommand{\tran}{\frown\hspace{-2.2ex}|\hspace{1.9ex}}
\newcommand{\lialg}{\sL}
\newcommand{\la}[2]{\Lambda_{#1#2}}
\newcommand{\ad}{\mathrm{ad}}
\newcommand{\kad}{\ad^{\#}}
\newcommand{\wl}[1]{\vphantom{X}_{#1}{\G}}
\newcommand{\te}{\tilde{e}}
\newcommand{\notka}[1]{}
\def\tgr{{\bf t}}
\def\sgr{{\bf s}}
\def\fgr{{\bf f}}
\def\rhogr{{\boldsymbol \rho}}
\newcommand{\sA}{\mbox{$\mathcal A$}}
\newcommand{\sT}{\mbox{$\mathcal T$}}
\newcommand{\sB}{\mbox{$\mathcal B$}}
\newcommand{\sF}{\mbox{$\mathcal F$}}
\newcommand{\sO}{\mbox{$\mathcal O$}}
\newcommand{\sD}{\mbox{$\mathcal D$}}
\newcommand{\sS}{\mbox{$\mathcal S$}}
\newcommand{\sY}{\mbox{$\mathcal Y$}}
\newcommand{\sL}{\mbox{$\mathcal L$}}
\newcommand{\hY}{\mbox{$\hat{Y}$}}
\newcommand{\hS}{\mbox{$\hat{S}$}}
\newcommand{\hX}{\mbox{$\hat{X}$}}
\newcommand{\dif}{differential }
\newcommand{\gru}{groupoid }
\newcommand{\grus}{groupoids }
\newcommand{\ti}{\tilde}
\newcommand{\halden}{half density }
\newcommand{\haldens}{half densities }
\renewcommand{\top}{topological }
\newcommand{\Setrel}{\mbox{\rm SetRel}}

\newcommand{\cstardwa}{\mbox{$C^*_r(\Gamma\times\Gamma)$}}
\newcommand{\xl}{X^L}
\newcommand{\xr}{X^R}
\newcommand{\kot}{\Pi}
\newcommand{\Af}{Aff}
\begin{document}
\title{On the Poisson structures related to $\kappa$-Poincar\'e Group.}
\author{Piotr Stachura}
%\abstract{tralalalaalal}
\address{Faculty of Applied Informatics and Mathematics, Warsaw University of Life Sciences-SGGW,
ul Nowoursynowska 166, 02-787 Warszawa, Poland,  
e-mail: stachura@fuw.edu.pl}
%\date{}
\begin{abstract} It is shown that the Poisson structure related to $\kappa$-Poincar\'e group is dual to certain Lie algebroid structure, the related Poisson structure on 
the (affine) Minkowski space is described in a geometric way. 
\end{abstract}
\maketitle
%%%%%%%%%%%%%%%%%%%%%%%%%%%%%%%%%%%%%%%%%%%%%%%%%%%%%%%%%%%%%%%%%%%%%%%%%%%%%%%%
%%%%%%%%%%%%%%%%%%%%%%%%%%%%%%%%%%%%%%%%%%%%%%%%%%%%%%%%%%%%%%%%%%%%%%%%%%%%%%%
\section{Introduction}
It is, more or less, ``common knowledge'' that the quantum $\kappa$-Poincar\'e  Group \cite{kappa1} exists on a $C^*$-level, it is given by some bicrossproduct construction 
\cite{Majid}, \cite{VV} and it's  a quantization of certain Poisson-Lie structure \cite{SZ-P}. 
Despite these beliefs, no precise and explicit formulae (e.g. for coproduct of generators) are known to the author. This note is a by-product of work on the $C^*$-version of
the $\kappa$-Poincar\'e. It consists of two parts.  In the first one, it is shown that really the Poisson structure presented in \cite{SZ-P} 
is dual to certain  Lie algebroid structure;  this Lie algebroid is the Lie algebroid of a groupoid.  
The $C^*$-algebra of this  groupoid should be the $C^*$-algebra of the quantum $\kappa$-Poincar\'e Group 
(it turns out we are in situation described in \cite{PS-ax}).  I tried to underline geometric and structural aspects of the construction. Such a formulation is necessary to
study whether and in what sense $\kappa$-Poincar\'e group is a quantization of the Poincar\'e Group.
The second part  describes the Poisson version of the  ``$\kappa$-Minkowski'' (affine) space and its relation to the Poisson structure on the Poincar\'e Group. 
These are again rather simple observations (essentially this part  is almost contained in \cite{SZ-PH}). 
In the second part, too, I tried to clarify geometric picture and  present results in a coordinate-free form.

{\bf Notation for orthogonal Lie algebras.}
$(V,\eta)$ stands for a  real, finite dimensional vector space with a bilinear, symmetric  and nondegenerate  form $\eta$.
{\em An orthonormal basis} is  a basis $(v_\alpha)$ in $V$ such that
$\displaystyle \eta(v_\alpha, v_\beta)=\eta(v_\alpha,v_\alpha)\delta_{\alpha\beta}\,,\,\,\,|\eta(v_\alpha,v_\alpha)|=1$.
For a vector $v$ with $|\eta(v,v)|=1$ we  write $sgn(v)$ for $\eta(v, v)$.
By $\eta$ we  denote also the isomorphism $V\rightarrow V^*$ given by $< \eta(x), y >:=\eta(x,y)$. Using  this notation,  for any orthonormal basis $(v_\alpha)$ 
and any $x,y \in V$:
\begin{equation}\label{on-basis}\notka{on-basis}
I=\sum_\alpha sgn(v_\alpha)v_\alpha\mt\eta(v_\alpha)\,,\,\,%v=\sum_\alpha sgn(v_\alpha) \eta(v_\alpha,v) v_\alpha\,,\,\,
\ \ \eta(x,y)=\sum_\alpha sgn(v_\alpha) \eta(x,v_\alpha) \eta(v_\alpha, y)
\end{equation}
A subspace generated by vectors $v_1,\dots\,, v_k$ is  denoted by $<v_1,\dots\,,v_k>$ or $span\{v_1,\dots\,,v_k\}$; for a subset $S\subset V$ the symbol $S^\perp$  denotes 
{\em the orthogonal complement} of $S$, if $S=\{v\}$ we  write $v^\perp$ instead of $\{v\}^\perp$; the symbol $S^0\subset V^*$ stands for {\em the anihilator} of $S$.

For vectors $x,y\in V$ let $\Lambda_{xy}:=x\mt\eta(y)-y\mt\eta(x)$; for a basis $(v_\alpha)$ in $V$ we  write $\Lambda_{\alpha\beta}$ instead of  $\Lambda_{v_\alpha,v_\beta}$. 
Operators $\Lambda_{xy}$ satisfy:
\begin{equation}\label{lambda-komut} \notka{lambda-komut} 
\,[\Lambda_{xy},\Lambda_{zt}]=\eta(x,t)\Lambda_{yz}+\eta(y,z)\Lambda_{xt}-\eta(x,z)\Lambda_{yt}-\eta(y,t)\Lambda_{xz}
\end{equation}
and $so(\eta)=span\{\Lambda_{xy}:x,y\in V\}$.
If $W\subset V$ is a subspace then $\Lambda_W:=span\{\Lambda_{xy}: x,y \in W\}$ is a Lie subalgebra of $so(\eta)$; for a {\em null vector} $f\in V$ and a subspace $W\subset V$ 
the subspace $\Lambda_{Wf}:= span\{\Lambda_{wf}: w\in W\}$ is also a subalgebra;  
notice   that for  $g\in O(\eta)$ we have  $\Lambda_{gx,gy}=g \Lambda_{xy} g^{-1}=: \ad(g) (\Lambda_{xy})$

\noindent
We will use a bilinear, nondegenerate form $k:so(\eta)\times so(\eta)\rightarrow \R$ defined by: 
\begin{equation}
\label{k-def}\notka{k-def}
k(\Lambda_{xy},\Lambda_{zt}):=\eta(x,t)\eta(y,z)-\eta(x,z)\eta(y,t)
\end{equation}
It is easy to see that for $g\in O(\eta)$: $\ad(g)\in O(k)$ i.e.
$$k(g\Lambda_{xy}g^{-1},g\Lambda_{zt}g^{-1})=k(\Lambda_{xy},\Lambda_{zt})\,,\,g\in O(\eta)$$
(of course $k$ is proportional to the Killing form). By $\kad$ we denote the coadjoint representation
of $O(\eta)$ on $so(\eta)^*$: $\kad(g):=\ad(g^{-1})^*$.
If $k$ is the isomorphism $so(\eta)\rightarrow so(\eta)^*$ defined by the form $k$ then 
$$\kad(g) k(X)=k(\ad(g) X) \,\,,\,X\in so(\eta)$$
Let us also define  a bilinear form $\tilde{k}$ on $so(\eta)^*$ by: 
\begin{equation}\label{tildek}\notka{tildek}
\tilde{k}(\varphi,\psi):=k(k^{-1}(\varphi),k^{-1}(\psi))\,,\,\varphi,\psi\in so(\eta)^*
\end{equation}
so $\tilde{k}(\varphi,\psi)=<\varphi, k^{-1}(\psi)>$;  
again it is clear that if $g\in O(\eta)$ then $\kad(g)\in O(\tilde{k})$, and
$$\tilde{k}(\kad(g) k(X),k(Y))= k(\ad(g) X, Y)\,\,,\,X,Y\in so(\eta)$$

We will also need one-parameter groups; they are given by formulae:
\begin{equation}\label{expy}\notka{expy}
\begin{split}
\exp(\Lambda_{u f})& =I+\Lambda_{u f}-\frac{\eta(u,u)}{2} f\mt\eta(f)\,\,,\ \ \ \ \eta(u,f)=0=\eta(f,f)\\
\exp(\nu \Lambda_{st})& = I-P_{st} +\cosh(\nu) P_{st} + \sinh(\nu)\Lambda_{st}\,\,,\ \ \eta(s,s)=-1=-\eta(t,t)\,,\,\eta(s,t)=0\,,\,\nu\in\R \\
\exp(\nu \Lambda_{xy})& = I-P_{xy} +\cos(\nu) P_{xy} + \sin(\nu)\Lambda_{xy}\,\,,\ \ \eta(x,x)=\eta(y,y)=\pm 1\,,\,\eta(x,y)=0\,,\,\nu\in\R 
\end{split}
\end{equation}
where $P_{vw}$ denotes the orthogonal projection onto $<v,w>$. 
If vectors $v,w$ are orthogonal and  $|\eta(v,v)|=|\eta(w,w)|=1$, then:
\begin{equation}P_{vw} = sgn(v) v\mt\eta(v)+sgn(w) w\mt \eta(w)=-sgn(v) sgn(w) \Lambda_{vw}^2\end{equation}
%%%%%%%%
{\bf Poincar\'e Group.} 
Let $(V,\eta)$ be  a vector Minkowski space (signature of $\eta$ is $(+,-,\dots,-)$).
For a vector $v\in V$ with $\eta(v,v)\neq 0$ let $R_v$ denote the reflection across the hyperplane $v^\perp$, i.e. 
$R_v = I-\frac{2}{\eta(v,v)} v\mt\eta(v)$. The full orthogonal group $O(\eta)$ has four connected components:
$SO_0(\eta)$ -- the connected component of identity; $R_t SO_0(\eta)\,,\,\eta(t,t)>0$ -- the component containing time reflection;  
$R_s SO_0(\eta)\,,\,\eta(s,s)<0$ -- the component containing space reflection and $SO_1(\eta):=R_t R_s SO_0(\eta)\,,\,\eta(t,t)>0, \eta(s,s)<0$ -- 
the component reversing time and space orientation (but keeping the space-time orientation intact). 
In this paper the {\em  Poincar\'e Group} $P(\eta)$ will mean  the semidirect product $V\rtimes O(\eta)$ and the {\em restricted Poincar\'e Group} $P_0(\eta)$ is  $V\rtimes SO_0(\eta)$.
Elements $(w,g)$ of $P(\eta)$ act on $V$ by affine mappings: $(w,g)(v):=w+gv$ and the group law is just the composition of these mappings:
$(w,g)(u,h)=(w+gu, gh)$. Since $P(\eta)$ depends only on dimension $n$ of $V$ it will be also denoted by $P(n)$; also $O(\eta)$ will be denoted by $O(1,n-1)$.

\section{Poisson-Poincar\'e Group.} 
The particular Lie-Poisson structure on  Poincar\'e Group we are interested in was defined in \cite{SZ-P}; 
it is dual to a certain Lie algebroid structure.  The construction is as follows.

Let  $(V,\eta)$ be a  vector Minkowski space of dimension $n+2\,,\,n>1$ and  $G:=SO_0(\eta)$. 
Our Poisson-Poincar\'e Group will be realized as a subgroup of the semidirect product $\gotg^*\rtimes G$:
\begin{equation}
\label{semi-direct}\notka{semi-direct}
(\varphi,g)(\psi, h):=(\varphi+\kad(g) \psi, g h)
\end{equation}
%%%%%%%%%%%%
Notice that if  $H\subset G$ is a subgroup with a Lie algebra $\gothh\subset \gotg$ then $\gothh^0\times H$ is  a subgroup of $\gotg^*\rtimes G$.
%($\gothh^0$ is the anihilator of $\gothh$ in $\gotg^*$).

Let us choose a (spacelike) vector   $\sgr\in V$ with $\eta(\sgr,\sgr)=-1$ and define a subalgebra:
\begin{equation}\label{a-def}\notka{a-def}
\gota:=span\{\Lambda_{xy}\,,\,x,y\in \sgr^\perp\}=\{Y\in so(\eta): Y \sgr =0 \}
\end{equation}
It is straightforward to see that:
$\displaystyle \gota^\perp=span\{\Lambda_{x\sgr}\,,\,x\in \sgr^\perp\}$  and
$k(\Lambda_{x\sgr},\Lambda_{y\sgr})=\eta(x,y)\,,\,x,y\in \sgr^\perp$ i.e. $(\gota^\perp, k)$ is {\em an $n+1$ dimensional vector Minkowski space} and the same is true for 
$(\gota^0,\tilde{k})$.

Let   $\tilde{A}$ be the connected subgroup of $G$ with Lie algebra $\gota$: % and let $A$ be the  normalizer of $\tilde{A}$. 
$\tilde{A}=\{g\in G: g\sgr=\sgr\}\simeq SO_0(1,n)$ (i.e $\tilde{A}$ is the proper, orthochronous Lorentz group). 
{\em Therefore the subgroup $\gota^0\times \tilde{A}$ is $P_0(n+1)$}; this way we have identified $P_0(n+1)$ as a subgroup of the semidirect product $\gotg^*\rtimes G$.
For reasons which are related to the ``quantum version'' of our Poisson-Poincar\'e group, we will also consider the  
normalizer of $\tilde{A}$ in $G$  which will be denoted by $A$. It is easy to see that
\begin{equation}\label{nor}\notka{nor}
A:=\{ g\in G: g\sgr=d(g) \sgr\,,\, d(g)=\pm 1\}=\tilde{A}\cup \exp(\pi\Lambda_{u\sgr}) \tilde{A}= \tilde{A}\cup  \tilde{A} \exp(\pi\Lambda_{u\sgr})
\end{equation}
for any spacelike, normalized  vector $u\in \sgr^\perp$.

% {\bf Tutaj opis normalizatora}

% \begin{lem} For  spacelike vectors $u, v \in <s>^\perp$:
% $$\exp(\pi\Lambda_{u\sgr}) \tilde{A}=\exp(\pi\Lambda_{v\sgr}) \tilde{A}=A$$
% \end{lem}
Let us now compute the action of $\exp(\pi\Lambda_{u\sgr})$ on $\gota^0$:
$$\kad(\exp(\pi\Lambda_{u\sgr}))k(\Lambda_{x\sgr})=k(\exp(\pi\Lambda_{u\sgr})\Lambda_{x\sgr} \exp(-\pi\Lambda_{u\sgr}))\,,\ {\rm and}$$
$$\exp(\pi\Lambda_{u\sgr})\Lambda_{x\sgr} \exp(-\pi\Lambda_{u\sgr})=-\Lambda_{x\sgr}-2 \eta(x,u)\Lambda_{u\sgr}\,,\,\eta(u,u)=-1\,,\,x, u \in \sgr^\perp$$
{\em In this way what exactly is $\gota^0\times A$ depends on the dimension of $V$:} for $n+1$ -- even this is $P_0(n+1)$ extended by time reflection; for $n+1$ -- odd this is
$P_0(n+1)$ extended by space and time reflection. 

The Lie-Poisson structure on $\gota^0\times A$ depends on a choice of a timelike vector $\tgr\in V$ or, equivalently, on a splitting of 
$\gota^0$ into  into {\em space}: $span\{k(\Lambda_{u \sgr})\,,\,u\in <\sgr,\tgr>^\perp\}$ and {\em time}: $<k(\Lambda_{\tgr \sgr})>$. So let us choose a (timelike) vector 
$\tgr\in \sgr^\perp\,,\,\eta(\tgr,\tgr)=1$; denote $\fgr:=\tgr-\sgr$ and let us define subalgebras:
\begin{equation}
\begin{split}\label{bc-def}\notka{bc-def}
\gotc& :=span\{\Lambda_{x\fgr}: x\in  \sgr^\perp\}=span\{\Lambda_{y\fgr}: y\in \tgr^\perp\}\\
\gotb& :=span\{\Lambda_{xy}\,,\,x,y\in \tgr ^\perp\}=\{Y\in so(\eta) : Y \tgr  =0\}
\end{split}
\end{equation}

The Lie algebra $so(\eta)$ can be decomposed as (direct sums of vector spaces):
\begin{equation} so(\eta)=\gotc\oplus \gotb=\gotc\oplus \gota \end{equation}

Let   $B, C$ be connected subgroups of $G$ with Lie algebras $\gotb,\gotc$ respectively; then  $B=\{g\in G: g \tgr=\tgr\}\simeq SO(n+1)$.
Denote  $U:=<\sgr,\tgr>^\perp\subset V$; then $(U,-\eta)$ is an $n$ dimensional (vector) Euclidean space. The subalgebra $\gotc$ can be decomposed further as
$$\gotc=\Lambda_{U\fgr}\ms <\Lambda_{\tgr \sgr}>,$$ 
where by (\ref{lambda-komut})  the first summand is an abelian ideal (in $\gotc$).\\
Using (\ref{expy}) we obtain: 
$$\exp(\nu \Lambda_{\tgr\sgr}) \exp(\Lambda_{u\fgr}) \exp(-\nu\Lambda_{\tgr\sgr})=\exp(\Lambda_{(e^\nu u) \fgr})\,,\,\, u\in U\,,\,\nu\in\R$$ 
Therefore $C=\{\exp(\Lambda_{u\fgr})\exp(\nu \Lambda_{\tgr\sgr}) : u\in U ,\,\nu\in\R\}$ and is isomorphic to the semidirect product $U\rtimes \R$ with group operation:
\begin{equation}\label{C-oper}\notka{C-oper}
(u,\mu) (v, \nu):=(u+e^\mu v, \mu+\nu)\,\,,\,u,v\in U\,,\,\mu,\nu\in \R
\end{equation}

The group  $C$ is the $AN$ group in the Iwasawa decomposition $SO_0(1,n+1)=SO(n+1) (AN)$ i.e there is the equality  $G=BC$.

The open set $\Gamma:=AC\cap CA$ carries two differential groupoid structures over $A$ and $C$ \cite{SZ-DG}.
The groupoid ``responsible'' for our Lie-Poisson structure is the groupoid $\G_A:\G\rightrightarrows A$.
Namely, the bundle $(TA)^0\subset T^*G$ is dual to the Lie algebroid $\lialg(\Gamma_A)$, which we realize as a vectors tangent in points of  $A$ to 
{\em left fibers} with bracket coming from {\em left invariant vector fields}. In this way  $(TA)^0$ carries the canonical Poisson structure; on the other hand
via right translation we can identify  $(TA)^0$ with $\gota^0\times A$ i.e. with the Poincar\'e group; it turns out this is a Poisson structure described in \cite{SZ-P}.
%
% Applying the cotangent lift to the groupoid $\G_A$ we obatin
% symplectic groupoid $T^*\G$ over $(TA)^0$. In this situation $(TA)^0$ is a Poisson manifolds and the Poisson bracket is dual 
% to the Lie algebroid structure on $\lialg(\Gamma_A).$
Let us compute Poisson brackets explicitely.

{\bf Algebroid structure.}
The  map $A\times \gotc\ni (a,p)\mapsto a p\in T_a\G_A$ is a global trivialization of $\lialg(\G_A)$.\\
For $p\in \gotc$ let $\xl_p$    be the  left invariant vector field on $\G_A$ defined by:
\begin{equation} \xl_p(a):=a p
\end{equation}
These vector fields satisfy:
\begin{equation}\label{bracket}\notka{bracket}\,
[\xl_p,\xl_q]=\xl_{[p,q]}\,,\,\,\,p,q\in\gotc
\end{equation}
and the anchor map $ \kot^L: \lialg(\G_A) \rightarrow TA$  is given by 
\begin{equation}\label{anchor}\notka{anchor}
\kot^L(\xl_p)(a)=P_{\gota}(\ad(a) p) a,
\end{equation} 
where $P_{\gota}$ is the projection onto $\gota$ corresponding to the decomposition $\gotg=\gotc\oplus\gota$. Short computations give:
\begin{equation}\label{anchor-1}\notka{anchor-1}
P_{\gota} \ad(a) \Lambda_{x\fgr}=\ad(a)\Lambda_{x\tgr}-d(a)\Lambda_{(ax)\tgr}\,,\,x\in \sgr^\perp \,,\,a\in A\,,\,a\sgr=:d(a) \sgr
\end{equation}

{\bf The Poisson structure.}
Sections of $\lialg(\G_A)$  define linear functions on $(TA)^0$, if $X$ is a section  of $\lialg(\G_A)$, 
the corresponding function will be denoted by $\tilde{X}$. Explicit form of this function for $\xl_p$ is:
$$\widetilde{\xl_p}(\varphi,a)= <\varphi a, \xl_p(a)>= <\varphi a, a p>=<\varphi, \ad(a) p>=\tilde{k}(\varphi,\kad(a) k(p))\,\,,\,\varphi\in \gota^0\,,\,\,p\in\gotc$$
(in this formula $(TA)^0\simeq \gota^0\times A$ via right translations).
%
%\noindent
The Poisson structure on $(TA)^0$ is defined by  the brackets:
\begin{equation}\label{poisson-bracket}\notka{poisson-bracket}
\{\widetilde{X_1}, \widetilde{X_2}\}=\widetilde{[X_1, X_2]}\,\,,\,\{\tilde{X}, \pi^*(f_1)\}=\pi^*(\kot^L(X) f_1)\,,\,\, \{\pi^*(f_1), \pi^*(f_2)\}=0, 
\end{equation}
where $f_1,f_2$ are smooth functions on $A$, $\pi: T^*G\rightarrow G$ is the cotangent bundle projection and $\pi^*$ denotes the pullback of functions.

Our Poincar\'e Group was identified with $\gota^0\times A\simeq (TA)^0$ (via right translations).
For $\varphi,\psi  \in \gota^0$ let us define smooth functions $\tilde{k}_\varphi, \tilde{k}_{\varphi\psi}$ on $\gota^0\times A$:
\begin{equation}
\begin{split}\
\tilde{k}_\varphi(\rho,a)& :=\tilde{k}(\varphi,\rho)\\
\tilde{k}_{\varphi\psi}(\rho,a)& :=\tilde{k}(\varphi, \kad(a) \psi)
\end{split}\end{equation}
Any Poisson structure on $\gota^0\times A$ is determined by brackets:
$$\{\tilde{k}_\varphi, \tilde{k}_\psi\}\,,\,\{\tilde{k}_\varphi, \tilde{k}_{\psi\rho}\}\,,\,\{\tilde{k}_{\varphi\lambda}, \tilde{k}_{\psi\rho}\}\,\,,
\,\,\varphi,\psi,\rho, \lambda\in \gota^0$$
For the Poisson structure given by (\ref{poisson-bracket}) we immediately get: 
\begin{equation}\label{nawiastryw}\notka{nawiastryw}
\{\tilde{k}_{\varphi\lambda}, \tilde{k}_{\psi\rho}\}=0.
\end{equation}
Let us now compute remaining brackets and compare them with the ones presented  in \cite{SZ-P}. To this end we will relate functions $\tilde{k}_\psi$ and $\widetilde{\xl_p}$.
\begin{lem}\label{lema1}\notka{lema1}
Let $(\rho_\alpha)$ be an orthonormal basis in $\gota^0$ and assume that elements $c_\alpha\in\gotc$ satisfy \\
$\tilde{k}(\psi,\kad(a)\rho_\alpha)=<\psi, \ad(a) c_\alpha>$ for any $\psi\in \gota^0$ and any $a\in A$.
Then:
$$\tilde{k}_\varphi=\sum_\alpha sgn(\rho_\alpha) \tilde{k}_{\varphi \rho_\alpha} \widetilde{\xl_{c_\alpha}}$$
\end{lem}
\noindent {\em Proof:} Indeed, using (\ref{on-basis}) let us compute:
\begin{equation*}
\begin{split}
\tilde{k}_\varphi(\psi, a)& =\tilde{k}(\varphi,\psi)=\tilde{k}(\kad(a^{-1}) \varphi,\kad(a^{-1}) \psi)=
\sum_\alpha sgn(\rho_\alpha)\tilde{k}(\kad(a^{-1}) \varphi,\rho_\alpha)\tilde{k}(\rho_\alpha,\kad(a^{-1}) \psi)=\\
 &= \sum_\alpha sgn(\rho_\alpha)\tilde{k}( \varphi,\kad(a)\rho_\alpha)\tilde{k}( \psi, \kad(a)\rho_\alpha)=
\sum_\alpha sgn(\rho_\alpha)\tilde{k}_{\varphi\rho_\alpha}(\psi, a) \tilde{k}( \psi, \kad(a)\rho_\alpha)=\\
& = \sum_\alpha sgn(\rho_\alpha)\tilde{k}_{\varphi\rho_\alpha}(\psi, a) <\psi, \ad(a) c_\alpha>=
\sum_\alpha sgn(\rho_\alpha)\tilde{k}_{\varphi\rho_\alpha}(\psi, a) \widetilde{\xl_{c_\alpha}}(\psi,a)=\\
& =\left(\sum_\alpha sgn(\rho_\alpha)\tilde{k}_{\varphi\rho_\alpha} \widetilde{\xl_{c_\alpha}}\right) (\psi,a)
\end{split}
\end{equation*}
\dowl

\noindent
Let $(v_\alpha)$ be an orthonormal basis in $\sgr^\perp$, then $(\rho_\alpha):=(k(\Lambda_{v_\alpha \sgr}))$ is an orthonormal  basis in $\gota^0$. 
Straightforward computations prove that elements
$c_\alpha:=-\Lambda_{v_\alpha\fgr}$ satisfy condition stated in the lemma above.

\noindent
Now, using (\ref{poisson-bracket}) and the decomposition above, we have:
\begin{equation}\label{nawiasy0}\notka{nawiasyO}
\begin{split}
\{\tilde{k}_\varphi, \tilde{k}_\psi\} & =
\sum_{\alpha\,\beta} sgn(\rho_\alpha) sgn(\rho_\beta) \{ \tilde{k}_{\varphi\rho_\alpha} \widetilde{\xl_{c_\alpha}}\,,\,\tilde{k}_{\psi\rho_\beta} \widetilde{\xl_{c_\beta}} \}=\\
& =\sum_{\alpha\,\beta} sgn(\rho_\alpha) sgn(\rho_\beta) \left[ \{ \widetilde{\xl_{c_\alpha}}, \tilde{k}_{\psi\rho_\beta} \}\,\tilde{k}_{\varphi\rho_\alpha} \widetilde{\xl_{c_\beta}}
- \{ \widetilde{\xl_{c_\beta}}, \tilde{k}_{\varphi\rho_\alpha} \}\,\tilde{k}_{\psi\rho_\beta} \widetilde{\xl_{c_\alpha}}+
\widetilde{\xl_{[c_\alpha,c_\beta]}} \tilde{k}_{\varphi\rho_\alpha} \tilde{k}_{\psi\rho_\beta}\right]=\\
&=:\framebox{{\rm I}}+\framebox{{\rm II}}+\framebox{{\rm III}}
\end{split}
\end{equation}
To end our computations we need formula for $\kot^L(\xl_p)(\tilde{k}_{\varphi\psi})$ for  $p:=\Lambda_{x\fgr}\in\gotc\,,\,x\in \sgr^\perp$ (note that the same symbol 
$\tilde{k}_{\varphi\psi}$ is used for function on $\gota^0\times A$ and on $A$). By (\ref{anchor}) and (\ref{anchor-1}):
\begin{equation}\label{Z-def}\notka{Z-def}
\kot^L(\xl_p)(a)= Z a\,\,,\,{\rm where}\,\,\, Z:=\ad(a)\Lambda_{x\tgr}-d(a)\Lambda_{(ax)\tgr}\,,\,a\in A\,,\,a\sgr=:d(a) \sgr
\end{equation}
and 
\begin{equation*}\begin{split}
\kot^L(\xl_p)(\tilde{k}_{\varphi\psi})(a)& =\frac{d}{d t}\Big\vert_{t=0} \tilde{k}_{\varphi\psi}(\exp(Z t) a)=\frac{d}{d t}\Big\vert_{t=0} \tilde{k}(\varphi,\kad(\exp(Z t) a) \psi)=\\
& = \frac{d}{d t}\Big\vert_{t=0} k(\Lambda_{v\sgr},\ad(\exp(Z t)) \ad( a) \Lambda_{w\sgr})=\frac{d}{d t}\Big\vert_{t=0} k(\ad(\exp(-Z t))\Lambda_{v\sgr}, \ad( a) \Lambda_{w\sgr}),
\end{split}\end{equation*}
where we put $\varphi=k(\Lambda_{v\sgr}), \psi=k(\Lambda_{w\sgr})$ for $v,w\in\sgr^\perp$.
We have the equality  $\displaystyle \ad(\exp(-Z t))\Lambda_{v\sgr}=\Lambda(\exp(-Z t)v, \exp(-Z t)\sgr)$, where for a while we use $\Lambda(x,y)$ for $\Lambda_{xy}$.
Since we are ineterested only in derivative in $t=0$ we can replace $\exp(-Z t)$ by $I-Zt$ and get:
\begin{equation}\label{der}\notka{der}
\frac{d}{d t}\Big\vert_{t=0} k(\ad(\exp(-Z t))\Lambda_{v\sgr}, \ad( a) \Lambda_{w\sgr})=k(-\Lambda_{(Zv)\sgr}-\Lambda_{v(Z\sgr)}, \ad(a) \Lambda_{w\sgr})
\end{equation}
By (\ref{Z-def}) $Z\sgr=0$ and :
$$Zv=\left[\Lambda_{(ax)( a\tgr)}-d(a)\Lambda_{(ax)\tgr}\right] v=\left[\eta(a\tgr,v)-d(a)\eta(\tgr,v)\right] (a x) - \eta(a x, v) (a\tgr)+ d(a) \eta( ax,v) \tgr.$$
Therefore
\begin{equation*}\begin{split}
\Lambda_{(Zv)\sgr}+\Lambda_{v(Z\sgr)}& =\Lambda_{(Zv)\sgr}=
\left[\eta(a\tgr,v)-d(a)\eta(\tgr,v)\right]\Lambda_{(ax)\sgr}-\eta(ax,v)\Lambda_{(a\tgr)\sgr}+d(a)\eta(ax,v)\Lambda_{\tgr\sgr}\\
&= \left[\eta(a\tgr,v)-d(a)\eta(\tgr,v)\right]d(a) \Lambda_{(ax)(a\sgr)}-\eta(ax,v) d(a)\Lambda_{(a\tgr)(a\sgr)}+d(a)\eta(ax,v)\Lambda_{\tgr\sgr}=\\
&= \left[\eta(a\tgr,v)-d(a)\eta(\tgr,v)\right]d(a) \ad(a)\Lambda_{x\sgr}-\eta(ax,v) d(a)\ad(a) \Lambda_{\tgr\sgr}+d(a)\eta(ax,v)\Lambda_{\tgr\sgr}
\end{split}\end{equation*}
So (\ref{der}) is equal to
\begin{equation}\label{der-cont}\notka{der-cont}
-\left[\eta(a\tgr,v)-d(a)\eta(\tgr,v)\right]d(a) \, k(\Lambda_{x\sgr},\Lambda_{w\sgr})+ \eta(ax,v) d(a)\, k(\Lambda_{\tgr\sgr},\Lambda_{w\sgr})-
d(a)\eta(ax,v)\,k(\Lambda_{\tgr\sgr},\ad(a) \Lambda_{w\sgr})
\end{equation}
Let us define $\rho:=k(\Lambda_{x\sgr})$ and $\rhogr:=k(\Lambda_{\tgr\sgr})$, then we have (recall that $\varphi=k(\Lambda_{v\sgr}), \psi=k(\Lambda_{w\sgr})$):
$$k(\Lambda_{x\sgr},\Lambda_{w\sgr})=\tilde{k}(\rho,\psi)\,,\,\,\,k(\Lambda_{\tgr\sgr},\Lambda_{w\sgr})=\tilde{k}(\rhogr,\varphi)\,,\,\,
k(\Lambda_{\tgr\sgr},\ad(a) \Lambda_{w\sgr})=\tilde{k}(\rhogr,\kad(a)\psi)=\tilde{k}_{\rhogr\psi}(a),$$
$$d(a)\eta(ax,v)=d(a)k(\Lambda_{(ax)\sgr},\Lambda_{v\sgr})=k(\ad(a)\Lambda_{x\sgr},\Lambda_{v\sgr})=\tilde{k}(\kad(a)\rho,\varphi)=\tilde{k}_{\varphi\rho}(a)$$
$$d(a)\eta(a\tgr,v)=d(a)k(\Lambda_{(a\tgr)\sgr},\Lambda_{v\sgr})=\tilde{k}_{\varphi\rhogr}(a)\,,\,\,\,\,\eta(\tgr,v)=\tilde{k}(\rhogr,\varphi)$$
and (\ref{der-cont}) is equal to:
\begin{equation*}
\begin{split}
& \left[\tilde{k}(\rhogr,\varphi)- \tilde{k}_{\varphi\rhogr}(a)\right] \tilde{k}(\rho,\psi)+ \tilde{k}_{\varphi\rho}(a)\left[\tilde{k}(\rhogr,\varphi)-\tilde{k}_{\rhogr\psi}(a)\right]=\\
& = \left\{\tilde{k}(\rho,\psi)\left[\tilde{k}(\rhogr,\varphi) I - \tilde{k}_{\varphi\rhogr}\right]+\tilde{k}_{\varphi\rho}\left[\tilde{k}(\rhogr,\varphi)I -
\tilde{k}_{\rhogr\psi}\right]\right\}(a)
\end{split}\end{equation*}
In this way we finally get:
\begin{equation}\label{kotnafun}\notka{kotnafun}
\kot^L(\xl_p)(\tilde{k}_{\varphi\psi})=\tilde{k}(\rho,\psi)\left[\tilde{k}(\rhogr,\varphi) I - \tilde{k}_{\varphi\rhogr}\right]+\tilde{k}_{\varphi\rho}\left[\tilde{k}(\rhogr,\varphi)I -
\tilde{k}_{\rhogr\psi}\right],
\end{equation}
where   $p:=\Lambda_{x\fgr}\,,\,x\in\sgr^\perp\,\,, \,\rho:=k(\Lambda_{x\sgr})$ and $\rhogr:=k(\Lambda_{\tgr\sgr})$.

Now we return to computations of (\ref{nawiasy0}). Choose an orthonormal basis $(e_\alpha)$ in $\sgr^\perp$ with $e_0:=\tgr$. Then we have orthonormal basis 
$\rho_\alpha:=k(\Lambda_{e_\alpha\sgr})=:k(\Lambda_{\alpha\sgr})$ in $\gota^0$ with $\rho_0=\rhogr$ and corresponding elements $c_\alpha:=-\Lambda_{e_\alpha\fgr} =: -\Lambda_{\alpha\fgr}$. Using (\ref{kotnafun}) we obtain:
\begin{equation*}\begin{split}
\{ \widetilde{\xl_{c_\alpha}}, \tilde{k}_{\psi\rho_\beta} \}& =\kot^L(\xl_{c_\alpha})(\tilde{k}_{\psi\rho_\beta})=\tilde{k}(-k(\Lambda_{\alpha\sgr}),k(\Lambda_{\beta\sgr}))\left[\tilde{k}(\rhogr,\psi)I-
\tilde{k}_{\psi\rhogr}\right]-\tilde{k}_{\psi\rho_\alpha}\left[ \tilde{k}(\rhogr,\psi)I  -\tilde{k}_{\rhogr\rho_\beta} \right]=\\
&=-\tilde{k}(\rho_\alpha,\rho_\beta)\left[\tilde{k}(\rhogr,\psi)I-
\tilde{k}_{\psi\rhogr}\right]-\tilde{k}_{\psi\rho_\alpha}\left[ \tilde{k}(\rhogr,\psi)I  -\tilde{k}_{\rhogr\rho_\beta} \right]=\\
&= -sgn(\rho_\alpha)\delta_{\alpha\beta}\left[\tilde{k}(\rhogr,\psi)I-
\tilde{k}_{\psi\rhogr}\right]-\tilde{k}_{\psi\rho_\alpha}\left[ \tilde{k}(\rhogr,\psi)I  -\tilde{k}_{\rhogr\rho_\beta} \right]
\end{split}
\end{equation*}
In this way the first term in the sum (\ref{nawiasy0})  is equal to:
\begin{equation*}
\begin{split}
\framebox{{\rm I}}=
\sum_{\alpha\,\beta} sgn(\rho_\alpha) sgn(\rho_\beta) \{ \widetilde{\xl_{c_\alpha}}, \tilde{k}_{\psi\rho_\beta} \}\,\tilde{k}_{\varphi\rho_\alpha} \widetilde{\xl_{c_\beta}} & =
-\tilde{k}_\varphi\left(\tilde{k}(\rhogr,\psi) I -\tilde{k}_{\psi\rhogr}\right) - \sum_{\alpha} sgn(\rho_\alpha)\tilde{k}_{\varphi\rho_\alpha}\tilde{k}_{\psi\rho_\alpha}\widetilde{\xl_{c_0}}+\\
& + \tilde{k}_\rhogr \sum_{\alpha} sgn(\rho_\alpha)\tilde{k}_{\varphi\rho_\alpha}\tilde{k}_{\psi\rho_\alpha}
\end{split}\end{equation*}
The second  term in (\ref{nawiasy0})  we get by intechanging in $\framebox{{\rm I}}$ $\alpha$ with $\beta$,  $\varphi$ with $\psi$ and changing the sign:
\begin{equation*}
\framebox{{\rm II}}=
\tilde{k}_\psi\left(\tilde{k}(\rhogr,\varphi) I -\tilde{k}_{\varphi\rhogr}\right) + \sum_{\beta} sgn(\rho_\beta)\tilde{k}_{\psi\rho_\beta}\tilde{k}_{\varphi\rho_\beta}\widetilde{\xl_{c_0}}-
\tilde{k}_\rhogr \sum_{\beta} sgn(\rho_\beta)\tilde{k}_{\psi\rho_\beta}\tilde{k}_{\varphi\rho_\beta}
\end{equation*}
and their sum is 
\begin{equation}
\framebox{{\rm I}}+ \framebox{{\rm II}}=\tilde{k}_\psi\left(\tilde{k}(\rhogr,\varphi) I -\tilde{k}_{\varphi\rhogr}\right)-\tilde{k}_\varphi\left(\tilde{k}(\rhogr,\psi) I -\tilde{k}_{\psi\rhogr}\right)
\end{equation}
It remains to compute $\framebox{{\rm III}}$.
$$[c_\alpha, c_\beta]=[-\Lambda_{\alpha\fgr},-\Lambda_{\beta\fgr}]=\eta(f,v_\beta)\Lambda_{\alpha\fgr}-\eta(f,v_\alpha) \Lambda_{\beta\fgr}=\delta_{0\alpha}c_\beta-\delta_{0\beta}c_\alpha$$
therefore
$$\widetilde{\xl_{[c_\alpha,c_\beta]}}=\delta_{0\alpha}\widetilde{\xl_{c_\beta}}-\delta_{0\beta}\widetilde{\xl_{c_\alpha}}$$
and
$$\framebox{{\rm III}}=\sum_{\alpha\,\beta} sgn(\rho_\alpha) sgn(\rho_\beta)\tilde{k}_{\varphi\rho_\alpha} \tilde{k}_{\psi\rho_\beta}
\left[\delta_{0\alpha}\widetilde{\xl_{c_\beta}}-\delta_{0\beta}\widetilde{\xl_{c_\alpha}}\right]=\tilde{k}_\psi \tilde{k}_{\varphi\rhogr}-\tilde{k}_\varphi \tilde{k}_{\psi\rhogr}$$
Finally:
\begin{equation}\label{nawias1}\notka{nawias1}
\{\tilde{k}_\varphi, \tilde{k}_\psi\}=\framebox{{\rm I}}+ \framebox{{\rm II}}+\framebox{{\rm III}}=\tilde{k}(\rhogr,\varphi) \tilde{k}_\psi-\tilde{k}(\rhogr,\psi) \tilde{k}_\varphi
\end{equation}
In the similar way, using lemma \ref{lema1} and formulae (\ref{nawiastryw}) and (\ref{kotnafun}) we obtain:
\begin{equation*}
\{\tilde{k}_{\lambda}, \tilde{k}_{\varphi \psi}\}   =\tilde{k}_{\lambda \psi}( \tilde{k}_{\varphi \rhogr} - \tilde{k}(\rhogr,\varphi)I ) + 
\tilde{k}(\lambda,\varphi)(\tilde{k}_{\rhogr \psi}- \tilde{k}(\rhogr ,\psi)I)
\end{equation*}

\noindent Now we have all the brackets:
\begin{equation}\label{brackets}\notka{brackets}
\begin{split}
\{\tilde{k}_{\varphi}, \tilde{k}_{\psi}\}& =\tilde{k}(\rhogr,\varphi) \tilde{k}_{\psi}- \tilde{k}(\rhogr,\psi) \tilde{k}_{\varphi}, \\
\{\tilde{k}_{\lambda}, \tilde{k}_{\varphi \psi}\}& = \tilde{k}_{\lambda \psi}( \tilde{k}_{\varphi \rhogr} - \tilde{k}(\rhogr,\varphi)I ) + 
\tilde{k}(\lambda,\varphi)(\tilde{k}_{\rhogr \psi}- \tilde{k}(\rhogr ,\psi)I), \\  
\{\tilde{k}_{\varphi\lambda}, \tilde{k}_{\psi\rho}\}& =0 \,\,\,{\rm for}\,\,\varphi, \lambda, \psi, \rho\in \gota^0\,\,{\rm and}\, \,\rhogr:=k(\Lambda_{\tgr \sgr}).
\end{split}
\end{equation}

The Poincar\'e group in \cite{SZ-P} was identified with matricies $g=\left( \begin{array}{cc} \Lambda, v\\0,1\end{array}\right) $, where $\Lambda$ is a Lorentz matrix of 
dimension $n+1$ and $v\in\R^{n+1}$. Poisson brackets for matrix elements of $g$ are given by:
\begin{equation}\label{rel-zak}\notka{rel-zak}
\begin{split}
\{\Lambda_{\mu\nu}, v_\beta\} & = h \left[(\Lambda_{\mu 0}-\delta_{\mu 0})\Lambda_{\beta\nu}+\eta_{\mu\beta}(\Lambda_{0\nu}-\delta_{0\nu})\right],\\
\{v_\alpha,v_\beta\}& = h( v_\alpha \delta_{\beta 0}-v_\beta\delta_{\alpha 0}),\\
\{\Lambda_{\mu\nu}, \Lambda_{\alpha\beta}\} & = 0,
\end{split}\end{equation}
where $\eta_{\alpha\beta}:=diag(1,-1,\dots,-1)$ and $h$ is a real parameter ({\bf Note:} here $\Lambda_{\alpha\beta}$ are {\em matrix elements not operators}).
To compare the  brackets, let us choose an orthonormal basis $(\rho_\alpha) \in \gota^0$ with $\rho_0=\rhogr$. 
We have 
$$\tilde{k}(\rho_\alpha,\rho_\beta)=diag(1,-1,\dots,-1)=\eta_{\alpha\beta}\,,\,\,  
v_\alpha=sgn(\rho_\alpha) \tilde{k}_{\rho_\alpha}\,\,{\rm and}\,\, \Lambda_{\alpha\beta}=sgn(\rho_\alpha)\tilde{k}_{\rho_\alpha\rho_\beta}.$$

Short computations show that brackets (\ref{brackets}) coincide with (\ref{rel-zak}) for $h=-1$.

\section{Poisson Minkowski space}

Let $(V,\eta)$ be a real, $n$-dimensional ($n>2$) vector space with a symmetric, bilinear, nondegenerate form $\eta$. For a basis $(v_\alpha)$ of $V$  let
$\eta_{\alpha\beta}:=\eta(v_\alpha,v_\beta)$ be the corresponding matrix of $\eta$ and $\eta^{\alpha\beta}$ stands for the inverse matrix. 
{\em Note that despite the title of the section, $(V,\eta)$ needn't to be a (vector) Minkowski space.}
In this section $G$ denotes any subgroup of $O(\eta)$ containing $SO_0(\eta)$ and $IG:=V\rtimes G$ is the semi-direct product. The Lie algebra of $IG$ is 
$iso(\eta):=V\times so(\eta)$ and the bracket is:
\begin{equation}
\label{in-bracket} \notka{in-bracket}\,
[(v,A), (w,B)]=(A w - B v, [A,B])
\end{equation}
The Poisson bracket for $\kappa$-Poincar\'e in \cite{SZ-P} is an example of a more general situation \cite{SZ-PPG}.
For a vector  $v\in V$ let us define 
\begin{equation}\label{def-b} b_v:=\sum \eta^{jk}e_j\wedge \Lambda_{v,e_k}\in iso(\eta)\wedge iso(\eta),\end{equation}
where $(e_k)$ is any basis in $V$.
Direct computation proves that, for $u,v\in V$, elements $b_v, b_u$ satisfy: 
\begin{equation}\label{b} [b_v, b_u]=-\eta(v,u) \Omega,\end{equation} 
%{\bf sprawdzone!!!} $[A ,v]=A v$ !!!!
where  $\Omega:=\sum \eta^{jk}\eta^{mn}e_j\wedge e_m \wedge \Lambda_{e_k,e_n}$ is the canonical %$ISO(\eta)$ 
invariant element in  $iso(\eta)\wedge iso(\eta) \wedge iso(\eta)$, and 
\begin{equation}\label{schouten}[a\wedge b, c\wedge d]:=
a\wedge [b,c]\wedge d -a\wedge [b,d]\wedge c- b\wedge [a,c]\wedge d+  b\wedge [a,d]\wedge c\end{equation}
is the (algebraic) Schouten bracket.
Therefore $b_v$ defines a Poisson-Lie  structure $\widehat{\Pi}_v$ on $IG$ by: %$ISO(\eta)$ by:
\begin{equation}\label{def-pi1} \widehat{\Pi}_v(g)=b_v g- g b_v\end{equation}
The structure in \cite{SZ-P} is of this type for $v$ being a timelike vector. Moreover, it is easy to see that
\begin{equation}\label{bv-wedge} 
[b_v, x\wedge u]=2 u\wedge x\wedge v\,,\,{\rm for\, any}\,\,x,u\in V
\end{equation} 
so we can replace $b_v$ in (\ref{def-pi1}) by $b_v+x\wedge v$ 
and  we obtain another Poisson-Lie structure on $IG$ which will be denoted by $\widehat{\Pi}_{v,x}$.
%%%%%%%%%
The adjoint representation of $IG$ on  $iso(\eta)$ is given by:
$$\ad_{(w,A)}(v,X)=(w+A v-A X A^{-1}w, A X A^{-1})\,,\,w,v\in V\,,\,A\in O(\eta)\,,\,X\in so(\eta);$$
by the same symbol we will denote this representation canonicaly extended to $iso(\eta)\wedge iso(\eta)$.\\
Straightforward computations give:
\begin{equation}\label{ad-b} \ad_{(w,A)}(b_v)=w\wedge Av + b_{Av}\,,\,\qquad\ad_{(w,A)}(x\wedge v)=Ax\wedge Av\end{equation}

Let $(M,V,\eta)$ be an affine space modeled on $(V,\eta)$. Let $\Af(G)$ be the group of those affine isometries of $M$ that have 
$G$ as their linear part.
Any point $m\in M$ defines the isomorphism $\phi_m: IG \rightarrow \Af(G)$ given by:
$$\phi_m(w,A)(m+v):=m+w+Av\,\,,\,\,v\in V$$
For two points $m,n\in M$ we have:  $\phi_m^{-1}\phi_n=Ad_{n-m}: IG\ni g\mapsto (n-m)\, g\, (n-m)^{-1}\in IG$ -- the inner automorphism given by $n-m\in V$.
%%%%%%%%%%%%
In this way for a point $m\in M$ and a vector $v\in V$ we have the Poisson-Lie structure on $\Af(G)$ defined by:
\begin{equation}\label{def-P}\Pi_{m,v}:=\phi_m(\widehat{\Pi}_v)\end{equation}
\begin{prop} Let $\Pi_{m,v}$ be the Poisson structure defined in (\ref{def-P}). Then:
\begin{itemize}
\item $\Pi_{m,\lambda v}=\lambda \Pi_{m,v}\,,\,\,\Pi_{m+\lambda v, v}=\Pi_{m,v}$ i.e.  the bivector $\Pi_{m,v}$  
depends only  on a parametrized line $l:=\{m+t v\,,\,t\in \R\}$; we will  write $\Pi_l$ for this Poisson structure.
\item Let $l, k$ be two parametrized lines then $\Pi_l=\Pi_k$ iff $l=k$.
\item If  $dim(V)>3$ then $\Pi_l$ and $\Pi_k$ are compatible iff $l$ and $k$ intersect or are parallel; \\
      for  $dim(V)=3$: if $G\subset SO(\eta)$ then  any two structures $\Pi_l$ and $\Pi_k$ are compatible;  otherwise the statement  is as for $dim(V)>3$.
\end{itemize}
\end{prop}
{\em Proof:} The equality $\Pi_{m,\lambda v}=\lambda \Pi_{m,v}$ is obvious. Let $m,n\in M$, $v, u\in V$ and  $x:=n-m\in V$. We can transfer $\Pi_{n,u}$ to  $IG$ by $\phi_m^{-1}$ and get
$\phi_m^{-1}\phi_n(b_u)=\ad_x(b_u)=x\wedge u+ b_u$ by (\ref{ad-b}). Taking $n:=m+\lambda v$ we get the second equality.

Let lines  $l, k$ be  given by $(m,v)$ and  $(n,u)$ respectively. Then $l\neq k$ means that $v\neq u$ or if $v=u$ then $x:=n-m\neq 0$ and $x,v$ are linearly independent.
Using the definition (\ref{def-b}) and the formula above it is easy to prove the second statement.

Poisson structures  $\Pi_l$ and $\Pi_k$ are compatible iff 
$\widehat{\Pi}_v+\widehat{\Pi}_{u,x}$ ($x:=n-m$) is a Poisson structure on $IG$, i.e. the Schouten bracket 
$[\widehat{\Pi}_v+\widehat{\Pi}_{u,x},\widehat{\Pi}_v+\widehat{\Pi}_{u,x}]=0$ and 
(since $\widehat{\Pi}_v$ and $\widehat{\Pi}_{u,x}$ are Poisson) this is equivalent to 
$[\widehat{\Pi}_v, \widehat{\Pi}_{u,x}]=0$. By (\ref{def-pi1}) this, in turn,  is equivalent to $[b_v, x\wedge u+b_u]$ being $IG$ invariant (with respect to adjoint action).
Using (\ref{b}) and (\ref{bv-wedge}) we get that $x\wedge u\wedge v$ must be $G$ invariant. 
Clearly this element is $0$ for intersecting or parallel lines  $l$ and $k$. For $dim(V)>3$ $G$,  invariance of $x\wedge u\wedge v$ forces it to be 0, i.e. 
lines $l$ and $k$ intersect or are parallel;  for $dim(V)=3$ the element $x\wedge u\wedge v$ is invariant if $G$ preserves orientation.\dowl

A parametrized line $l:=\{m+t w\,,\,t\in \R\}$ defines also a bivector $\pi_l$ on $M$: 
\begin{equation}\label{pi-M}\pi_l(m+v):=v\wedge w,
\end{equation}
in the formula above we identify $TM$ with $M\times V$; it is  easy to see that really $\pi_l$ depends only on $l$ and not on the chosen point $m\in l$. 
\begin{prop}
\begin{itemize}
\item $\pi_l$ is a Poisson bivector on $M$.
\item $\pi_l$ and $\pi_k$ are compatible iff lines $l, k$ intersect or are parallel.
\item The canonical action of $(\Af(G),\Pi_k)$ on $(M,\pi_l)$ is Poisson iff $l=k$.
\end{itemize}
\end{prop}
%%%%%%%
{\em Proof:} Let $l:=\{m+t w\,,\,t\in \R\}$ and define the vector field $\hat{V}^m$ by $\hat{V}^m(m+v):=v$; let $\hat{w}$ be the constant vector field:
$\hat{w}(m+v):=w$; with this notation we have: $\pi_l=\hat{V}^m\wedge \hat{w}$. If $\pi_k$ is defined by the line $k:=\{n+t u,\,t\in\R\}$ then
$$\pi_k(m+v)=\pi_k(n+(m-n)+v)=(m-n)\wedge u + v\wedge u=( \hat{x}\wedge\hat{u}+ \hat{V}^m\wedge \hat{u})(m+v), $$
where $x:=m-n$, i.e. $\pi_k=\hat{x}\wedge\hat{u}+ \hat{V}^m\wedge \hat{u}$. Let us compute:
%%%
$$[\pi_l,\pi_k]=[\hat{V}^m\wedge \hat{w},\hat{x}\wedge\hat{u}+ \hat{V}^m\wedge \hat{u} ]=
[\hat{V}^m\wedge \hat{w},\hat{x}\wedge\hat{u}] + [\hat{V}^m\wedge \hat{w},\hat{V}^m\wedge \hat{u} ]=$$
$$=-\hat{w}\wedge [\hat{V}^m,\hat{x}]\wedge \hat{u}+  \hat{w}\wedge [\hat{V}^m,\hat{u}]\wedge \hat{x}+ 
\hat{V}^m\wedge [\hat{w},\hat{V}^m]\wedge \hat{u}+ \hat{w}\wedge [\hat{V}^m,\hat{u}]\wedge \hat{V}^m$$
But for any constant vector field $\hat{y}$ we have: $[\hat{V}^m,\hat{y}]=-\hat{y}$, therefore:
$$[\pi_l,\pi_k]=2 \hat{w}\wedge \hat{x}\wedge \hat{u}.$$
In this way  $[\pi_l,\pi_l]=0$ and  $\pi_l+\pi_k$ is Poisson iff $w\wedge x\wedge u=0$. Now first and the second statement are clear.

Let the lines $l,k$ be defined by $(m,w)$ and $(n,u)$, respectively; let $\psi_n:V\ni v\mapsto n+v\in M$. Using $\phi_n$ and $\psi_n$ we can transfer problem to the 
action on $(IG, \hat{\Pi}_u)$ on $(V,\hat{\pi}_l)$, where $\hat{\Pi}_u$ is defined by (\ref{def-pi1}) and $\psi_n(\hat{\pi}_l)=\pi_l$ i.e. $\hat{\pi}_l(v)=(x+v)\wedge w\,,\,x:=n-m$.
The action is 
$$IG\times V\ni (y,A;v)\mapsto y+Av\in V$$
This action is Poisson iff
\begin{equation}
\label{P-cond}
\hat{\pi}_l(gv)=\hat{g}\hat{\pi}_l(v)+ \hat{\Pi}_u(g)\hat{v}\,,\,\,g:=(y,A)\in IG, 
\end{equation}
where $\hat{g}$ is (the extension of) the mapping $V\ni v\mapsto gv\in V$ and $\hat{v}$ (the extension of)
$IG\ni g\mapsto gv\in V$. 

\noindent
We have:
$$\hat{\pi}_l(gv)=\hat{\pi}_l(y+Av)=(x+ y+ Av)\wedge w$$
$$\hat{g}(\hat{\pi}_l(v))=(Ax+Av)\wedge Aw$$
$$\hat{\Pi}_u(g)\hat{v}=(b_u g - g b_u)\hat{v}=(b_u)\widehat{gv}-\hat{g}(b_u\hat{v})$$
It is straightforward, that for $(\dot{x},\dot{A})\in T_eIG$ : 
$(\dot{x},\dot{A})\hat{z}=\dot{x}+\dot{A}z$; so 
$$(b_u)\hat{z}=\sum \eta^{jk}e_j\wedge (\Lambda_{u e_k} z)=\sum \eta^{jk}e_j\wedge (\eta(e_k,z) u - \eta(u,z)e_k)=z\wedge u$$
therefore
 $$(b_u)\widehat{gv}=(y+Av)\wedge u\,\,,\,\,\,(b_u)\widehat{v}=v\wedge u$$
and $$\hat{g}(b_u\hat{v})=\hat{g}(v\wedge u)=Av\wedge Au$$
In this way equality (\ref{P-cond}) reads:
$$(x+ y+ Av)\wedge w=(Ax+Av)\wedge Aw+(y+Av)\wedge u-Av\wedge Au\,\,,\,\,{\rm for\,any\,}\, y,v\in V\,,\,\,A\in G$$
If  $l=k$ i.e. $x=0, w=u$ this condition is fulfilled. On the other hand, 
setting $v=0, A=I$ we get (for any $y$) $y\wedge w=y\wedge u$, so $w=u$ and the equality reduces to
$$x\wedge w=Ax\wedge Aw\,\, {\rm for\,any\,}\, A\in G$$
Therefore $x\wedge w=0$ and  $l=k$.
\dowl


\begin{thebibliography}{}

\bibitem{kappa1} J Lukierski, A Nowicki, H Ruegg, {\em New quantum Poincaré algebra and κ-deformed field theory} Physics Letters B, {\bf 293} (1992), pp 344-352.
\bibitem{Majid} S Majid, H Ruegg, {\em  Bicrossproduct structure of κ-Poincaré group and non-commutative geometry}, Physics Letters B, {\bf 334} (1994), pp 348-354.
\bibitem{VV} S. Vaes, L. Vainerman, {\em  Extensions of locally compact quantum groups and the bicrossed product construction}, Advances in Mathematics {\bf 175} (1) (2003), 1--101 . 
\bibitem{SZ-P} S. Zakrzewski {\em Quantum Poincar\'{e} group related to the $\kappa$-Poincar\'{e} algebra.}  J. Phys. A Math. Gen. {\bf 27} (1994) 2075-2082.
\bibitem{PS-ax} P. Stachura {\em On the quantum ‘ ax + b ’ group}, J. Geom. Phys. {\bf 73} (2013), 125-149.
\bibitem{SZ-PPG} S. Zakrzewski {\em Poisson Structures on the Poincar\'e Group}, Comm. in Math. Phys. {\bf 185} (1997), pp 285-311.
\bibitem{SZ-PH} S. Zakrzewski, {\em Poisson homogeneous spaces}, in: J. Lukierski, Z. Popowicz, J. Sobczyk (eds.), Quantum groups (Karpacz, 1994), PWN, Warszawa, 1995,
pp 629–639.
\bibitem{SZ-DG} S. Zakrzewski, {\em Quantum and classical pseudogroups. II. Differential and symplectic pseudogroups}, Comm. Math. Phys. {\bf 134},  (1990), pp 71-395.
\end{thebibliography}
\end{document}